\documentclass[11pt]{amsart}

\usepackage{amssymb}
\usepackage{amsmath}
\usepackage{graphicx}
\usepackage{tikz}
\usepackage{appendix}
\usepackage{bibentry}

       \font\tenmsb=msbm10
       \font\sevenmsb=msbm7
       \font\fivemsb=msbm5
       \catcode`\@=11
       \ifx\amstexloaded@\relax\catcode`\@=\active
       \endinput\else\let\amstexloaded@\relax\fi
       \def\spaces@{\space\space\space\space\space}
       \def\spaces@@{\spaces@\spaces@\spaces@\spaces@\spaces@}
       \def\space@.{\futurelet\space@\relax}
       \space@. %
       \def\Err@#1{\errhelp\defaulthelp@\errmessage{AmS-TeX error: #1}}
       \def\relaxnext@{\let\next\relax}
       \def\accentfam@{7}
       
       \def\noaccents@{\def\accentfam@{0}}
       \def\Cal{\relaxnext@\ifmmode\let\next\Cal@\else
       \def\next{\Err@{Use \string\Cal\space only in math mode}}\fi\next}
       \def\Cal@#1{{\Cal@@{#1}}}
       \def\Cal@@#1{\noaccents@\fam\tw@#1}

       \def\Bbb{\relaxnext@\ifmmode\let\next\Bbb@\else
       \def\next{\Err@{Use \string\Bbb\space only in math mode}}\fi\next}
       
       \def\Bbb@#1{{\Bbb@@{#1}}}
       \def\Bbb@@#1{\noaccents@\fam\msbfam#1}
       \newfam\msbfam
       \textfont\msbfam=\tenmsb
       \scriptfont\msbfam=\sevenmsb
       \scriptscriptfont\msbfam=\fivemsb

\def\Z{{\Bbb Z}}

\newtheorem{thm}{Theorem}[section]

\newtheorem{lem}[thm]{Lemma}

\theoremstyle{remark}
\newtheorem{rem}{Remark}[section]

\theoremstyle{definition}

\newtheorem{definition}{Definition}[section]

\newenvironment{thmbis}[1]
  {\renewcommand{\thethm}{\ref{#1}$'$}%
   \addtocounter{thm}{-1}%
   \begin{thm}}
  {\end{thm}}

\@addtoreset{equation}{section}

\newcommand{\beq}{\begin{equation} }
\newcommand{\eeq}{\end{equation} }

\begin{document}
\setlength{\columnsep}{5pt}
\title{Uniform Localization is always Uniform}

\author{Rui Han}
\address{
University of California, Irvine, California
}
\email{rhan2@uci.edu}

\thanks{$^\dag$This work was partially supported in part by DMS-1401204.}

\date{}
\maketitle
\begin{abstract}
In this note we show that if a family of ergodic Schr$\ddot{o}$dinger operators on $l^2(\Z^\gamma)$ with continuous potentials have uniformly localized eigenfunctions then these eigenfunctions must be uniformly localized in a homogeneous sense.
\end{abstract}

\section{Introduction}
\noindent

Given a topological space $\Omega$, let $T_i : \Omega\rightarrow\Omega$ be commuting homeomorphisms, and let $\mu$ be an ergodic Borel measure on $\Omega$. Let $f : \Omega\rightarrow \mathbb{R}$ be continuous and define $V_\omega(n)=f(T^n\omega)$ for $n\in \Z^\gamma$, where $T^n=T_1^{n_1}...T_\gamma^{n_\gamma}$. Let $H_\omega$ be the operator on $l^2(\Z^\gamma)$,
$$(H_\omega u)(n)=\sum_{|m-n|=1}{u(m)+V_\omega(n) u(n)}.$$

The occurrence of pure point spectrum for the operators $\{H_\omega\}$ is called $phase\, stable$ if it holds for every $\omega\in\Omega$.

For a self-adjoint operator $H$ on $l^2(\Z^\gamma)$. We say that $H$ has uniformly localized engenfunctions ($ULE$), if $H$ has a complete set of orthonormal eigenfunctions $\{\phi_n\}_{n=1}^\infty$, and there are $\alpha>0$, $C>0$, such that
$$|\phi_n(m)|\leq Ce^{-\alpha|m-m_n|}$$
for all eigenfunctions $\phi_n$ and suitable $m_n$. It is known that $ULE$ has a close connection with phase stability of pure point spectrum. Actually in paper \cite{Jit97}, Jitomirskaya pointed out that instability of pure point spectrum implies absence of uniform localization. It is also shown in \cite{DRJLS96} that if $H_\omega$ has $ULE$ for $\omega$ in a set of positive $\mu$-measure, then $H_\omega$ has pure point spectrum for any $\omega\in supp(\mu)$. The proof of this statement mainly relies on the fact that $ULE$ implies uniform dynamical localization ($UDL$), which means if $H_{\omega}$ has $ULE$, then
$$|(\delta_l, e^{-itH_\omega}\delta_m)|\leq C_\omega e^{-\alpha_\omega|l-m|}$$
for some constants $\alpha_\omega$, and $C_\omega$ that depend on $\omega$. Recently in \cite{DG10} and \cite{DG11}, Damanik and Gan established $ULE$ for a certain model and then proved that for this model, actually $C_\omega$ and $\alpha_\omega$ can be chosen to be independent of $\omega$. In this note we will show that the latter property is always a corollary of $ULE$.

First, let us give a new definition.

\begin{definition}
$H_\omega$ has uniform or homogeneous $ULE$ in a set $S$ means that $H_\omega$ has $ULE$ for any $\omega$ in $S$ and
$$|\phi^\omega_n(m)|\leq Ce^{-\alpha|m-m^\omega_n|}$$
with constants $\alpha>0$ and $C>0$ which do not depend on $\omega$.
\end{definition}

Then the main theorems in the note can be stated as follows:

\begin{thm}\label{comm1}  If $H_{\omega}$ has $ULE$ for $\omega$ in a positive $\mu$-measure set, then $H_\omega$ has homogeneous $ULE$ in $supp(\mu)$.
\end{thm}

\begin{thmbis}{comm1} If $T$ is minimal, and $H_\omega$ has $ULE$ at a single $\omega$, then $H_\omega$ has homogenous $ULE$ in $\Omega$.
\end{thmbis}



\section{Proof of Theorem 1.1}

Let $\{\mathbf{U}_k\}$ be a family of transitions on $l^2(\mathbb{Z}^\gamma)$ defined by $(\mathbf{U}_k u)(m):=u(m-k)$. Clearly, if $\{\phi_n^\omega\}$ is a complete set of eigenfunctions of $H_\omega$, then $\{\mathbf{U}_k \phi_n^\omega\}$ is a complete set of eigenfunctions of $H_{T^k\omega}$. Also, if $H_\omega$ has $ULE$, which means $\exists\ \alpha_0 >0, C_0>0$ such that
$$|\phi_n^\omega(m)|\leq C_0e^{-\alpha_0|m-m_n^\omega|}$$
for all eigenfunctions $\phi_n^\omega$ and suitable $m_n^\omega$. Then $H_{T^k\omega}$ also has $ULE$. In fact if we let $m_n^{T^k\omega}=m_n^\omega+k$, then $$|\phi_n^{T^k\omega}(m)|=|(\mathbf{U}_k \phi_n^\omega)(m)|\leq C_0 e^{-\alpha_0|m-m_n^{T^k\omega}|}$$
for all eigenfunctions $\mathbf{U}_k \phi^\omega_n$, also notice that the constants $C_0$ and $\alpha_0$ are the same for $H_\omega$ and $H_{T^k\omega}$.

\begin{lem} If $H_\omega$ has $ULE$ for $\omega$ in a positive $\mu$-measure set $S$, then $H_\omega$ has $ULE$ for $a.e.$ $\omega\in\Omega$. \end{lem}

\noindent \emph{Proof}. $H_\omega$ has $ULE$ in $\bigcup\limits_{k\in \mathbb{Z}^\gamma}T^k S$, which is a transition invariant set, so $\mu(\Omega\backslash(\bigcup\limits_{k\in \mathbb{Z}^\gamma}T^k S))=0$.\qed

\begin{thm}  If $H_\omega$ has $ULE$ for $a.e.\ \omega\in\Omega$, then $\exists\ \alpha >0$ independent of $\omega$, such that
$$|\phi_n^\omega(m)|\leq C_\omega e^{-\alpha|m-m_n^\omega|}$$
for $a.e.\ \omega\in\Omega$ and all eigenfunctions $\phi_n^\omega$ with suitable $m_n^\omega$.
\end{thm}

\noindent \emph{Proof}. Let $\bigcup\limits_{j=1}^\infty\{\omega|\ |\phi_n^\omega(m)|\leq C_\omega e^{-\frac{1}{j}|m-m_n^\omega|},\ $
for all eigenfunctions
$\phi_n^\omega$ and suitable $m_n^\omega$\}$:=\bigcup\limits_{j=1}^{\infty} A_j$.
\
$A_j$ is translation invariant. Since $\mu(\bigcup\limits_{j=1}^{\infty} A_j)=1$, $\exists\ j_0$ such that $\mu(A_{j_0})=1$. \qed

\

\noindent Now, let's return to the proof of Theorem 1.1.\\

\noindent \emph{Proof of Theorem 1.1}.
By Theorem 2.2, and direct computation we have $|(\delta_l,e^{-itH_\omega}\delta_{m})|\leq C_\omega e^{-\alpha|l-m|}$ for a.e. $\omega$. Let
$$F(\omega)=\sup_{t\in \mathbb{Q},\ l,m\in \mathbb{Z}^\gamma}|(\delta_l, e^{-itH_\omega}\delta_m)|e^{\alpha|l-m|}.$$
Then $F(\omega)<\infty\ a.e.\ \omega$. It is also easy to see that $F(\omega)$ is measurable and translation invariant. Therefore, by ergodicity, $F(\omega)=C\ a.e.\ \omega$. Hence $|(\delta_l, e^{-itH_\omega}\delta_m)|\leq Ce^{-\alpha|l-m|}\ a.e.\ \omega$. Then on a dense set in $supp(\mu)$,
$$|(\delta_l, e^{-itH_\omega}\delta_m)|\leq Ce^{-\alpha|l-m|}.$$
By continuity, the inequality holds for any $\omega$ in $supp(\mu)$.

\noindent Then since $P^{\omega}_{\{E\}}=$s-$\lim_{T\rightarrow\infty}\frac{1}{2T}\int^{T}_{-T}e^{iEs}e^{-iH_\omega s}ds$, we have
$$|(\delta_l, P^\omega_{\{E\}}\delta_m)|\leq Ce^{-\alpha|l-m|}$$
for any $E\in \mathbb{R}$. Therefore if we choose $\tilde{m}_n^\omega$ so that $|\phi^\omega_n(\tilde{m}_n^\omega)|=\sup_m|\phi^\omega_n(m)|$, we get
$$|\phi^\omega_n(l)|^2\leq|\phi^\omega_n(l)||\phi^\omega_n(\tilde{m}_n)|\leq C e^{-\alpha|l-\tilde{m}^\omega_n|}.$$
\qed

\begin{rem}
For the proof of Theorem $1.1^\prime$, one need to realize that when $T$ is minimal, $F(\omega)$ being translation invariant implies that $F(\omega)$ is constant in a dense subset of $\Omega$.
\end{rem}

\section{Generalization}

In fact we can extend the result above to a more general case where $f(x)$ is allowed to have discontinuities.

\begin{definition}
We say $f$ has invariant continuity filter in $\Omega$ if at every $\omega\in\Omega$, there is a filter $F_\omega$, such that any $A_\omega\in F_\omega$ satisfies the following conditions:\\
$1.\ \mu(A_\omega\bigcap B(\omega,\delta))>0, $ for any $\delta>0$,\\
$2.\ \lim_{\omega_k\in A_\omega, \omega_k\rightarrow\omega}f(\omega_k)\rightarrow f(\omega),$\\
$3.\ T^n(A_\omega)\in F_{T^n\omega}, $ for any $n\in \mathbb{Z}^\gamma$.
\end{definition}

\noindent  Example. \ \
Let $\Omega=\mathbb{T}=\mathbb{R}/\mathbb{Z}$, and $\mu$ be the Lebesgue measure. For any $\theta\in\Omega$,  $T(\theta)=\theta+\alpha$ where $\alpha\notin \mathbb{Q}$ and $f(x)=\{x\}$. Define $(H_\theta u)(n)=u(n+1)+u(n-1)+f(T^n\theta)u(n)$. The reason why we are interested in this model is that $ULE$ has recently been shown for it in \cite{JKpre}. Obviously in this model $f$ is not continuous but it does have continuity invariant filter at every $\theta\in [0,1]$. In fact, the filter at $\theta$ is the set of all intervals with left endpoint $\theta$. Generally speaking, all the right or left continuous function defined on $\mathbb{R}$ with direction preserving $T$ has invariant continuity filter at every point.

\

Now we have the following theorem:
\begin{thm}\label{comm2} Assume f is bounded and has invariant continuity filter at every $\omega\in\Omega$, then if $H_{\omega}$ has $ULE$ in a positive $\mu$-measure set, $H_\omega$ has homogenous $ULE$ in $supp(\mu)$.
\end{thm}

As before, we also have:
\begin{thmbis}{comm2} Assume f is bounded and has invariant continuity filter at every $\omega\in\Omega$, then if $H_{\omega}$ has $ULE$ at a single $\omega$, $H_\omega$ has homogenous $ULE$ in $\Omega$.
\end{thmbis}

\noindent \emph{Proof}.
Notice that in the Proof of Theorem 1.1, we only use the continuity of $f$ in the last step, which means we still have
$$(\delta_l, e^{-itH_\omega}\delta_m)\leq Ce^{-\alpha|l-m|}$$
for $\omega\in\Omega_0$, where $\mu(\Omega_0)=1$.

\noindent Now consider any $\omega_0\notin\Omega_0$, we know $\mu(A_{\omega_0}\bigcap B(\omega_0, \frac{1}{k}))>0$, hence we can choose $\omega_k^{(0)}\in A_{\omega_0}\bigcap B(\omega_0, \frac{1}{k})$,\, $\omega_k^{(0)}\in\Omega_0$ and $|f(\omega_k^{(0)})-f(\omega_0)|<\frac{1}{2^k}$. Then $T^{m}\omega_k^{(0)}\in T^{m} A_{\omega_0}$, for any $m\in \mathbb{Z}^\gamma$. $T^m\omega_k^{(0)}\rightarrow T^m\omega_0$, therefore $f(T^m\omega_k^{(0)})\rightarrow f(T^m\omega_0)$. Hence we can choose a subsequence of $\{\omega_k^{(0)}\}$, say $\{\omega_k^{(m)}\}$, such that $|f(T^m\omega_k^{(m)})-f(T^m\omega_0)|<\frac{1}{2^k}$. Notice that by the diagonal argument, we can find a sequence $\{\omega_k\}$, satisfying $\omega_k\in\Omega_0$ and $|f(T^j\omega_k)-f(T^j\omega_0)|<\frac{1}{2^k}$ for any $j\in \mathbb{Z}^\gamma$ and $k\geq |j|$.

Now, let us show that for fixed $l, m, t$:
$$(\delta_l, e^{-itH_{\omega_k}}\delta_m)\rightarrow (\delta_l, e^{-itH_{\omega_0}}\delta_m).$$

\noindent Indeed: $|(\delta_l, e^{-itH_{\omega_k}}\delta_m)-(\delta_l, e^{-itH_{\omega_0}}\delta_m)|$\\
$=|(\delta_l, (e^{-it(H_{\omega_k}-H_{\omega_0})}-1)e^{-itH_{\omega_0}}\delta_m)|$\\
$=|(\delta_l, (e^{-it(H_{\omega_k}-H_{\omega_0})}-1)\sum_{r=-\infty}^{\infty}a_r\delta_r)|$\\
$\leq|a_l|(e^{|t||f(T^l\omega_k)-f(T^l\omega_0)|}-1)\rightarrow 0$, as $k\rightarrow\infty$.

\noindent Hence $|(\delta_l, e^{-itH_{\omega_0}}\delta_m)|\leq Ce^{-\alpha|l-m|}$.\qed

\section*{Acknowledgement}
I would like to thank Svetlana Jitomirskaya for her guidance and inspiring discussions on this subject.

\bibliographystyle{amsplain}

\end{document}